\input AHTOH-E.STY
\hfuzz8pt

\UDC{
512.542.1     
+ 512.543.1   
+ 512.552.18  
}


\title{%
Strange divisibility in groups and rings
}
\author{%
Anton A. Klyachko$^\sharp$
\quad
Anna A. Mkrtchyan$^\flat$
}
\address{
$^\sharp$\myAddressW
\\
$^\flat$University of Edinburgh,
School of Mathematics,
Room 5402, James Clerk Maxwell Building, King's Buildings,\\
Edinburgh, EH9 3JZ
\\
anna.mkr@gmail.com
}

\grantsFirst{\RFBR15-01-05823}

\abstract{%
We prove a general divisibility theorem that implies, e.g., that, in any
group, the number of generating pairs (as well as triples, etc.) is a
multiple of the order of the commutator subgroup. Another corollary says
that, in any associative ring, the number of Pythagorean triples (as well
as four-tuples, etc.) of invertible elements is a multiple of the order of
the multiplicative group.
}

\s 0.
Introduction

The starting point of our investigation is the following fact generalising
an old theorem of Solomon [Solo69].

\proclaim{Gordon--Rodriguez-Villegas Theorem} \rm [GRV12].
Let $F$ be a finitely generated group whose commutator subgroup is of
infinite index and let $G$ be an arbitrary group. Then the number of
homomorphisms $F\to G$ is divisible by the order of $G$.

This theorem is essentially about the number of solutions to systems
of coefficient-free equations in a group. In~[KM14], this result was
extended to equations with coefficients and even to arbitrary first-order
formulae in group language (with constants).

The main theorem of this paper has a claim to the title of ``the maximal"
generalisation of the Gordon--Rodriguez-Villegas Theorem (although such
a maximality can never be proven). The statement of
the Main Theorem can be found in the first section; a (quite elementary)
proof is in last section. Roughly speaking, the Main Theorem asserts that
the divisibility is retained when we take into account only some set of
homomorphisms provided the set is invariant with respect to some natural
operations on homomorphisms. One of the corollaries of the Main Theorem is
the unexpected fact mentioned in the abstract:
\disp{%
\hfuzz39pt
\sl
in any group $G$, the number of generating tuples
$(g_1,\dots,g_{\the\year})\in G^{\the\year}$ \(i.e. such tuples that
$G=\gp{g_1,\dots,g_{\the\year}}$\) is a multiple of the
order of the commutator subgroup of~$G$.
}
(Here, the number \the\year\ can be replaced by any integer;
see Section 2 for a more general fact and other group-theoretic
applications.) Surprisingly, this result seems to
be new, although many related facts on the divisibility of the M\"obius
function (which is related to the number of generating tuples via the Hall
formula [Hall36]) are known, see, e.g., [Bro00], [HI\"O89], [KT84], and
references therein.  We refer to~[Coll10] for yet other not widely known
but beautiful results about generating tuples.

The Main Theorem is an assertion about groups, but (paradoxically) it has
nontrivial ring-theoretic corollaries. In Section 3, we derive
a ring-theoretic analogue of the Gordon--Rodriguez-Villegas theorem (to be
more precise, an analogue of the generalisation of this theorem from
[KM14], which is about equations with coefficients). A particular
case of this theorem on equations over rings is the fact mentioned in the
abstract, or, e.g., the following higher-order assertion:
\disp{%
\sl
in any associative ring $R$ with unity, the number of tuples of
invertible elements
$(a,b,\dots,z)\in (R^*)^{26}$ such that
$
a^{\the\year}+b^{\the\year}+\dots+z^{\the\year}=0
$
is divisible by the order of the multiplicative group of this ring, i.e.
by~$|R^*|$.
}
(Here, the number \the\year\ can be replaced by any integer; see Section 3
for a more general fact.)

\smallskip
\noindent
{\bf Our notation}
is mainly standard. Note only that, if $k\in \Z$ and $x$ and $y$ are
elements of a group, then $x^y$, $x^{ky}$, and~$x^{-y}$ denote $y^{-1}xy$,
$y^{-1}x^ky$ and $y^{-1}x^{-1}y$, respectively.  The commutator subgroup
of a group $G$ is denoted $G'$.  If $X$ is a subset of a group, then
$|X|$, $\gp X$, $\nc X$, and $C(X)$ are the cardinality of~$X$, the
subgroup generated by~$X$, the normal closure of~$X$, and the centraliser
of~$X$, respectively. The index of a subgroup $H$ of a group $G$ is
denoted $|G:H|$.
The symbol~$N(H)$ denotes the normaliser of a subgroup $H$ (in a
group $G$). The free product of groups $A$ and $B$
is denoted $A*B$ and $F(x_1,\dots,x_n)$ is the free group with basis
$x_1,\dots,x_n$.
If $R$ is an
associative ring with unity, then $R^*$ denotes the group of invertible
elements of this ring.

Note also that, in almost all assertions about the divisibility
(e.g., in the Gordon--Rodriguez-Villegas Theorem), one need not
assume that the corresponding group is finite. The divisibility
can be understood in the sense of cardinal arithmetics: an
infinite cardinal is divisible by all nonzero smaller (and equal)
cardinals.
We really need finiteness assumption only in the Theorem on Monomorphisms
and Subgroups in Section 2 (see the corresponding remark there).

The authors thank E. B. Vinberg for a question that brought us to writing
this paper, Andrey~V.~Vasil'ev for a valuable remark (see
Section 2), and an anonymous referee for a lot of comments allowing us to
improve the text.

\s 1.
Main Theorem

A group $F$ equipped with an epimorphism $F\to\Z$ is called
\emph{indexed}. This epimorphism $F\to\Z$ is called \emph{degree} and
denoted $\deg$; thus, an integer $\deg f$ is assigned to each element $f$
of an indexed group $F$ in such a way that $F$ contains elements of all
integer degrees and ${\deg(fg)=\deg f+\deg g}$ for any $f,g\in F$.

Suppose that $\phi\:F\to G$ is a homomorphism from an indexed group $F$
to a group $G$ and $H$ is a subgroup of~$G$. We say that the subgroup
$$
H_\phi=\bigcap_{f\in F}H^{\phi(f)}\cap C(\{\phi(f)\;|\;\deg f=0\})
$$
is the \emph{$\phi$-core} of $H$. In other words, the $\phi$-core $H_\phi$
of~$H$ consists of such elements $h$ that $h^{\phi(f)}\in H$ for all $f$
and, moreover, $h^{\phi(f)}=h$ if $\deg f=0$.

\proclaim{Main Theorem}.
Let $H$ be a subgroup of a group $G$ and let $\Phi$ be a set of
homomorphisms from an indexed group~$F$ to $G$
with the following two properties.
\item{\rm I.}
$\Phi$ is invariant with respect to conjugation by elements of $H$:
$$
\qbox{if $h\in H$ and $\phi\in\Phi$, then the homomorphism
$\psi\:f\mapsto\phi(f)^h$ lies in $\Phi$.}
$$
\item{\rm II.}
For any $\phi\in\Phi$ and any $h$ from the $\phi$-core $H_\phi$
of $H$, the homomorphism $\psi$ defined by
$$
\psi(f)=
\cases{
\phi(f)& for all $f\in F$ of degree zero;
\cr
\phi(f)h& for
some element
$f\in F$ of degree one
\small(and, hence, for all degree-one elements)
\cr
}
$$
belongs to $\Phi$ too.
\enditem
Then $|\Phi|$ is divisible by $|H|$.

Note that the mapping $\psi$ from Condition I is a homomorphism for any
$h\in G$; and the formula for $\psi$ from Condition~II defines a
homomorphism for any $h\in C(\phi(\ker\deg))$
(see Lemma 0).
Conditions I and II only
require these homomorphisms to belong to $\Phi$ (under some additional
restrictions on $h$).

\Lemma 0.
Suppose that
$\phi\:F\to G$ is a homomorphism from an indexed group $F$
to a group $G$,
$f_1$ is a degree-one element of $F$,
and
$g\in G$.
Then
\item{\rm 1)}
if (and only if)
$g\in C(\phi(\ker\deg))$, then
there exists
a (unique) homomorphism $\psi\:F\to G$
such that $\psi(f)=\phi(f)$ for all $f$ of degree zero
and $\psi(f_1)=\phi(f_1)g$;
\item{\rm 2)}
if
$H$ is a subgroup of $G$
and $g\in H_\phi$, then
$\psi(f)H=\phi(f)H$ for all $f\in F$.

\Proof
Note that $F$ is a semidirect product
$F=\gp{f_1}_\infty\semitimes\ker\deg$.
This means that a mapping $\alpha\:\ker\deg\cup\{f_1\}\to G$
extends to a homomorphism if and only if its restriction to
$\ker\deg$ is a homomorphism and
$\alpha(f^{f_1})=\alpha(f)^{\alpha(f_1)}$ for all~$f\in\ker\deg$.

For all~$f\in\ker\deg$,
we have
$\psi(f^{f_1})=\phi(f^{f_1})=\phi(f)^{\phi(f_1)}$
and
$\psi(f)^{\psi(f_1)}=\phi(f)^{\phi(f_1)g}$.
This implies that
$\psi(f^{f_1})=\psi(f)^{\psi(f_1)}$ for all~$f\in\ker\deg$
if and only if $\phi(x)^g=\phi(x)$ for all~$x\in\ker\deg$.
This proves the first assertion.

To prove 2) note that any $f\in F$ has the form $f=f_1^kx$, where
$x\in\ker\deg$
and $k\in\Z$. So,
$$
\psi(f)H=\psi(f_1)^k\psi(x)H=\psi(f_1)^k\phi(x)H=(\phi(f_1)g)^k\phi(x)H
=\!=\!=
\phi(f_1)^k\phi(x)H=\phi(f_1^kx)H=\phi(f)H
$$
(where the equality $=\!=\!=$ is valid because
$\phi(F)$ normalises $H_\phi$ and $g\in H_\phi\subseteq H$).
This proves assertion~2).

\s 2.
Applications. Groups

First, note that the conditions of the Main Theorem are obviously
satisfied if $\Phi$ is the set of all homomorphisms $F\to G$ (and $H$ is
any subgroup of $G$, e.g., the entire group $G$). Therefore, the
Gordon--Rodriguez-Villegas theorem is the simplest special case of the
Main Theorem.

\Theorem on Equations over Groups \rm [KM14].
The number of solutions to a system of equations
\newline
${\{v_i(x_1,\dots,x_n)=1\}}$ over a group~$G$ \(where
$v_i(x_1,\dots,x_n)\in G*F(x_1,\dots,x_n)$\) is divisible by the order of
the centraliser of the set of coefficients if the rank of the matrix
consisting of the exponent-sums of $i$-th unknown in $j$-th equation is
less than the number of unknowns.

\Proof
Let $A\subseteq G$ be the subgroup generated by all coefficients of the
equations. Let $F$ be the quotient
group ${F=(A*F(x_1,\dots,x_n))/\nc{\{v_i\}}}$ of the free product
$A*F(x_1,\dots,x_n)$ of $A$ and the free group by the normal
subgroup~$\nc{\{v_i\}}$ generated by the left-hand sides of the equations.
Let $\Phi$ be the set of homomorphisms $F\to G$ that are identity on $A$.
(We assume that $A$ embeds into $F$ via the natural map $A\to F$, because
if this map is not injective, then there are no solutions and we have
nothing to prove.) Clearly, solutions to the system of equations are in a
natural one-to-one correspondence with the elements of $\Phi$.

The condition on the rank means that $F$ admits an epimorphism onto $\Z$
whose kernel contains $A$. Let $H$ be the centraliser of $A$ in $G$.
Clearly, the conditions of the Main Theorem are satisfied.  Indeed,
Condition I holds, because~$h$ centralises $A\subseteq G$ and, hence,
$\psi$ coincides with $\phi$ on $A\subset F$; Condition II holds, because
elements of $A\subset F$ are of degree zero and, hence, $\psi$ coincides
with $\phi$ on $A\subset F$ again.

\Theorem on Roots of Subgroups \rm [KM14].
The number of elements $g$ of a group $G$ such that $g^n\in H$ is
divisible by $|H|$ for any subgroup $H$ of $G$ and any integer $n$.%
\fn{\rm
In 2017, we learned that this fact was proven in [Iwa82].
}

Theorem on Roots of Subgroups is the simplest special case of the
following fact.

\Theorem on Homomorphisms and Subgroups \rm [KM14].
Let $H$ be a subgroup of a group $G$ and let $W$ be a subgroup (or subset)
of a finitely generated group $F$ whose commutator subgroup $F'$ is of
infinite index. Then the number of homomorphisms $\phi\:F\to G$ such that
$\phi(W)\subseteq H$ is divisible by $|H|$.

We shall prove a yet more general fact.

\Theorem on Homomorphisms and Double Cosets.
Let $H$ be a subgroup of a group $G$, let $W$ be a subset of a finitely
generated group $F$ whose commutator subgroup $F'$ is of infinite index,
and let $W\ni w\mapsto g_w\in G$ be an arbitrary map $W\to G$.  Then the
number of homomorphisms $\phi\:F\to G$ such that $\phi(w)\in Hg_wH$
for all $w\in W$ is divisible by $|H|$.

\Proof
Take some epimorphism $\deg\:F\to\Z$ (which exists because $F/F'$ is an
infinite finitely generated abelian group) and let $\Phi$ be the set of
all homomorphisms $\phi\:F\to G$ such that $\phi(w)\in Hg_wH$ for
all $w\in W$.  The conditions of the Main Theorem hold. For Condition I,
this is quite obvious. As for Condition~II, it suffices to note that the
formula for~$\psi$ implies the equality $\psi(f)H=\phi(f)H$ for all
$f\in F$ by Lemma 0.

\smallskip

The following theorem is an ``epimorphism analogue" of the
Gordon--Rodriguez-Villegas theorem.

\Theorem on Epimorphisms.
Let $F$ be a finitely generated group whose commutator subgroup is of
infinite index and let $G$ be an arbitrary group. Then the number of
surjective homomorphisms $F\to G$ is divisible by the order of the
commutator subgroup of~$G$.

\Proof
Take some epimorphism $\deg\:F\to\Z$, let $\Phi$ be the set of all
epimorphisms $F\to G$, and put $H=G'$. Let us verify that the conditions
of the Main Theorem are satisfied. For Condition I, this is obvious.

To verify Condition II, we have to show that, for any epimorphism
$\phi\:F\to G$ and any element~${h\in G'}$ centralising the subgroup
$\phi(\ker\deg)$, the homomorphism $\psi$ from Condition II is
surjective. Clearly, it is surjective modulo $G'$ (i.e. $\psi(F)G'=G$),
because $\psi$ equals $\phi$ modulo $G'$. It remains to show that each
element $g\in G'$ lies in $\psi(F)$. By the surjectivity of
$\phi$, we can find $f\in F$ such that $\phi(f)=g$; moreover, the
element~$f$ can be found in the commutator subgroup of $F$ (because, for
an epimorphism, the image of the commutator subgroup equals the commutator
subgroup of the image). But then $f\in\ker\deg$ and, therefore,
$\psi(f)=\phi(f)=g$ as required.

\Remark.
{
The number of surjective homomorphisms $F\to G$ is a multiple of
$|\Aut G|$, because $\Aut G$ acts faithfully on the set of epimorphisms
$F\to G$.  However, the Theorem on Epimorphisms does not follows
immediately from this observation, because, as was noted by A. V. Vasil'ev,
$$
\hbox{\sl
there exists a group $G$ such that $|\Aut G|$ is not divisible by $|G'|$.
}
$$
Examples of such groups are the groups ${3\cdot A_6}$ and ${3\cdot A_7}$
(see, e.g., [Wils09]) of orders ${{3\over2}\cdot6!=1080}$ and
${{3\over2}\cdot7!=7560}$; they coincides with their commutator subgroups
and have centres of order three; the central quotients are the alternating
group $A_6$ and $A_7$ while ${|\Aut\!(3\cdot A_6)|=2\cdot6!}$ and
$\Aut\!(3\cdot A_7)$ is the symmetric group of order~$7!$.
Actually, Savelii Skresanov and Dmitrii Churikov showed (using {\tt GAP})
that the smallest group $G$ such that $|G'|{\not|\;} |\Aut G|$ is of order
108.  }

\Corollary on Generating Tuples in Groups.
For each group $G$ and each positive integer $n$, the number of tuples
$(g_1,\dots,g_n)\in G^n$ of elements of $G$ generating $G$ \(i.e.  such
that $\gp{g_1,\dots,g_n}=G$\) is divisible by~$|G'|$.

\Proof
The generating $n$-tuples of elements of $G$ are in a natural one-to-one
correspondence with the epimorphisms from the free group of rank $n$ to
$G$. Therefore, the assertion follows immediately from the Theorem
on Epimorphisms.

\smallskip

Clearly, in the Theorem on Epimorphisms and even in the above corollary,
the divisibility by $|G'|$ cannot be strengthened to the divisibility
by $|G|$, because, in a prime-order group, the number of generating
$n$-tuples is $|G|^n-1$.

The following theorem generalises the Theorem on Epimorphisms and is
an analogue of the Theorem on Homomorphisms and Subgroups.

\Theorem on Epimorphisms and Subgroups.
Let $A$ be a subgroup of a group $G$ and let $W$ be a subgroup of a
finitely generated group $F$ whose commutator subgroup $F'$ is of infinite
index. Then the number of homomorphisms $\phi\:F\to G$ such that
$\phi(W)=A$ is divisible by $|A'|$.

\Proof
Take some epimorphism $\deg\:F\to\Z$ and put
$$
\Phi=\{homomorphisms\ \phi\:F\to G\ such\ that\ \phi(W)=A\}
\qqbox{and}
H=A'.
$$
Let us verify that the condition of the Main Theorem are satisfied. For
Condition I, it is obvious.

To verify Condition II, we have to show that, for any homomorphism
$\phi\:F\to G$ such that $\phi(W)=A$ and any element $h\in A'$
centralising the subgroup $\phi(\ker\deg)$, we have $\psi(W)=A$ for the
homomorphism $\psi$ from Condition II. The inclusion $\psi(W)\subseteq A$
certainly holds. To prove the inverse inclusion, note that $\psi(W)A'=A$
It remains to show that each element $a\in A'$
lies in $\psi(W)$. Since $\phi(W)=A$, we can find $w\in W$ such that
$\phi(w)=a$; clearly, such an
element $w$ can be found in the commutator subgroup of~$W$. But then
$w\in\ker\deg$ and, therefore, $\psi(w)=\phi(w)=a$ as required.

\medskip

A similar statement on injective homomorphisms also holds (for finite
groups $G$); moreover, the divisibility is much better in this case.

\Theorem on Monomorphisms and Subgroups.
Let $A$ be a subgroup of a group $G$ and let $W$ be a subgroup of a
finitely generated group $F$ such that $WF'$ is of infinite
index in $F$. Then $|N(A)|$ divides the following numbers:
\item{\rm a)}
the number of homomorphisms $\phi\:F\to G$ such that
the restriction of $\phi$ to $W$ is injective and
$\phi(W)\subseteq A$;
\item{\rm b)}
the number of homomorphisms $\phi\:F\to G$ such that
the restriction of $\phi$ to $W$ is injective and
$\phi(W)=A$.

\Proof
Let us prove a) (the proof of b) is quite similar).
Take an epimorphism $\deg\:F\to\Z$
such that $W\subseteq\ker\deg$
and put
$$
\Phi=\{homomorphisms\ \phi\:F\to G\ such\ that\ \phi(W)\subseteq A
\ and\ \phi|_W\ is\ injective\}
\qqbox{and}
H=N(A).
$$
Condition I of the main theorem is obviously satisfied.
Condition II is also satisfied, because
$W\subseteq\ker\deg$
and, hence, $\psi$ and
$\phi$ (from condition II) coincide on $W$.

\Remark.
The condition $|F:F'W|=\infty$ cannot be
replaced
by $|F:F'|=\infty$
in the last
theorem (even though we understand divisibility in the sense of cardinal
arithmetics). Indeed,
\item{\rm a)}
if $F=W=A=\Z$ and $G=\R$, then
the number of injective homomorphisms $F=W\to A$
is $\aleph_0$ which is not a multiple of $|N(A)|=|\R|=2^{\aleph_0}$;
\item{\rm b)}
if $F=W=G=A=\Z$, then
the number of bijective homomorphisms
is two which is not a multiple of $|N(A)|=|\Z|=\aleph_0$.

\s 3.
Applications. Rings

A \emph{generalised homogeneous} equation over an associative ring $R$
with the set of unknowns $X$ is a finite equation of the form
$$
\sum_i\prod_j c_{ij}x_{ij}^{k_{ij}}=0,
\qbox{where
\emph{coefficients} $c_{ij}\in R$,
\emph{unknowns} $x_{ij}\in X$,
and \emph{exponents} $k_{ij}\in\Z$},
$$
such that for some nonzero mapping $\deg\:X\to\Z$ the value
$\sum\limits_jk_{ij}\deg(x_{ij})$ does not depend on $i$ (i.e. the
``polynomial" in the left-hand side of the
equation is homogeneous with respect to
some nonzero assignment of degrees to variables%
\fn{%
A variable may have zero degree,
but at least one variable must have a nonzero degree.
}%
).
A system of equations is called generalised homogeneous if all its
equations are generalised homogeneous (possibly of different degrees) with
respect to the same function~$\deg\:X\to\Z$.

To test generalised homogeneity, one can use the following
simple algorithm.

\medskip

{\noindent\bf
Algorithm for Testing Generalized Homogeneity of a System}

\item{\bf1.}
For each equation $v=0$, construct the matrix $A_v$ with integer entries
$a_{ij}$ that are the degree of $i$-th monomial with respect to $j$-th
unknown (i.e. $a_{ij}$ is the exponent-sum of $j$-th unknown in $i$-th
monomial of $v$).

\item{\bf2.}
Subtract the first row of this matrix $A_v$ from each row of $A_v$.
Do it for all matrices $A_v$.

\item{\bf3.}
Combine the matrices $A_v'$ thus obtained (with zero first rows)
into one matrix:
$
A'=\pmatrix{
A_v'
\cr
A_w'
\cr
\vdots
\cr
}.
$

\item{\bf4.}
The system is generalised homogeneous if and only if the rank of $A'$ is
less than the number of unknowns.

\medskip

\noindent
For example, for the system of equations
$
\cases{
(xdy)^2-yx^2+xy^2cy^{-100}x=0
\cr
xy-yx=0
}
$
(where $c,d\in R$ are coefficients and $x,y$ are unknowns),
we obtain:
$$
A_u=\pmatrix{
2&2\cr
2&1\cr
2&-98\cr
},\
A_v=\pmatrix{
1&1\cr
1&1\cr
},
\quad
A_u'=\pmatrix{
0&0\cr
0&-1\cr
0&-100\cr
},\
A_v'=\pmatrix{
0&0\cr
0&0\cr
},
\quad
A'=\pmatrix{
0&0\cr
0&-1\cr
0&-100\cr
0&0\cr
0&0\cr
},
$$
$\rank A'=1$
and the system is generalised homogeneous.

\Proposition.
Any system of equations such that
$$
\sum_i
\Big((\hbox{the number of monomials in $i${\rm th} equation})-1\Big)
<(\hbox{the number of unknowns})
$$
is generalised homogeneous.

\Proof
The assertion follows immediately from the above algorithm, but we leave
the proof of correctness of this algorithm to readers as an exercise. (We
shall use neither this proposition, nor this algorithm in this paper.)

\medskip

The notion of a \emph{solution} to a system of equations is defined
naturally (if some exponents $k_{ij}$ are negative, then the corresponding
components of the solution must be invertible elements of the ring).

\Th on Equations over Rings.
Let $R$ be an associative ring with unity and let $G$ be a subgroup of the
multiplicative group of this ring. Then, for each generalised homogeneous
system of equations over~$R$ with $n$ unknowns, the number of solutions
lying in $G^n$ is divisible by the order of the intersection of $G$ and
the centraliser of the set of coefficients of the system.

\Proof
Let us apply the Main Theorem by letting $F$ be the free group
$F(x_1,\dots,x_n)$ and extending the mapping
$\deg\:\{x_1,\dots,x_n\}\to\Z$ (from the definition of generalised
homogeneous systems) to a homomorphism $F\to\Z$, which can be assumed to
be surjective, because it is nonzero. Let $\Phi$ be the set of
homomorphisms $\phi\:F\to G$ such that $(\phi(x_1),\dots,\phi(x_n))$ is a
solution to the system of equations, and let $H$ be the intersection of
$G$ and the centraliser of the set of coefficients of the system.

Let us verify the conditions of the Main Theorem. Condition I
holds obviously. To verify Condition~II, choose an element $t\in F$
of degree one and write each variable $x_i$ in the form
$x_i=t^{\deg x_i}y_i$, where ${y_i=t^{-\deg x_i}x_i}$ has zero degree.

Consider an equation $w(x_1,\dots,x_n)=0$ of the system and let us
rewrite it in the form~$v(t,y_1,\dots,y_n)=0$. By virtue of
homogeneity, all monomials in $v(t,y_1,\dots,y_n)$ have the same
degree $k$ with respect to $t$.

We have to show that, if $v(\phi(t),\phi(y_1),\dots,\phi(y_n))=0$ and
$h\in H_\phi$, then $v(\phi(t)h,\phi(y_1),\dots,\phi(y_n))=0$.
And this is
indeed so, because
$v(\phi(t)h,\phi(y_1),\dots,\phi(y_n))$ is a multiple of
$v(\phi(t),\phi(y_1),\dots,\phi(y_n))$ by the following lemma
(which should be applied to each monomial of $v$).

\Lemma 1.
Suppose that $M$ is a monoid, $b_i,a,h\in M$, elements $a$ and $h$ are
invertible,
and $a^{-s}ha^s$, where $s\in\Z$, commute with all $b_i$. Then,
for any expression of the form
$$
u(t)=b_0t^{n_1}b_1\dots t^{n_l}b_l,
\qbox{where $n_i\in\Z$,}
$$
we have
$
u(ah)=
\cases{
h^{a^{-1}}h^{a^{-2}}\dots h^{a^{-k}}u(a) & if $k=\sum n_i>0$
\cr
h^{-1}h^{-a}\dots h^{-a^{-1-k}}u(a) & if $k=\sum n_i<0$
\cr
u(a) & if $k=\sum n_i=0$.
\cr
}
$

\Proof
Using the commuting rules $a^ih^{a^j}=h^{a^{j-i}}a^i$ and
$b_ih^{a^j}=h^{a^j}b_i$, we bring all letters $h$ (and $h^{a^j}$) to the
left end of the word $u(ah)$ and obtain the required form. This
completes the proofs of Lemma 1 and the Theorem on Equations over Rings.

\medskip

\Example.
{\sl
The number of \emph{Pythagorean triples} of invertible elements
of an associative ring with unity, i.e. the number of invertible solutions
to the equation
$$
x^2+y^2=z^2
$$
is always divisible by the order of the multiplicative group of the ring.
}
\newline
Indeed, the equation is homogeneous and
we can take $G=R^*$. Moreover,
\newline
{\sl
the number of invertible solutions
to the equation
$$
ax^k+by^l+cz^m+dt^n+\dots=0
$$
is divisible by $|R^*|$ for any
$a,b,c,d,\dots,k,l,m,\dots\in\Z$,
}
\newline
because this equation is generalised homogeneous.


\s 4.
Proof of the main theorem

The argument is to some extent similar to that near the end of Section~3
of~[KM14]. To emphasise the similarity we use the same terms
as in ~[KM14]
(albeit their meaning is different).

The \emph{tail} of a homomorphism $\phi\in\Phi$ is the pair
$(\phi_0,\phi_H)$, where $\phi_0$ is the restriction of~$\phi$ to the
subgroup $\ker\deg\subset F$ and $\phi_H\: F\to\{gH\;;\;g\in G\}$ is
the mapping from $F$ to the set of left cosets of $H$ in $G$ that sends
an element $f\in F$ to the coset~$\phi(f)H$.

We say that two homomorphisms $\phi,\psi\in\Phi$ are \emph{similar} and
write $\phi\sim\psi$ if their tails are conjugate by an element of
$H$, i.e.
$$
\eqalign{
\phi\sim\psi\ \iff
\qqbox{there exists $h\in H$
such that}
&\hbox{$\psi(f)=h\phi(f)h^{-1}$}
\qqbox{for all $f\in F$ of degree zero and}
\cr
&\psi(f)H=h\phi(f)H
\!\qbox{for all $f\in F$.}
}
$$
Clearly, similarity is an equivalence relation on $\Phi$. The Main
Theorem is an immediate corollary of the following proposition.

\Proposition.
In $\Phi$, each class of similar homomorphisms consists of exactly $|H|$
elements. More precisely, for each~${\phi\in\Phi}$,
\item{\rm 1)}
the number of different tails of elements of~$\Phi$
similar to $\phi$ is $|H:H_\phi|$;
\item{\rm 2)}
for each homomorphism $\psi$ similar to $\phi$, the number of elements of
$\Phi$ with the same tail as $\psi$ is~$|H_\phi|$.

\Proof
To prove 1), note that the group $H$ acts by conjugation on the set of
tails of elements of~$\Phi$.  Indeed, if we conjugate the tail of a
homomorphism $\psi\in\Phi$ by an element $h\in H$, then we obtain the tail
of the homomorphism $f\mapsto\psi(f)^h$. This homomorphism lies in $\Phi$
by Condition I of the Main Theorem. The tails of homomorphisms similar
to~$\phi$ form the orbit of the tail of $\phi$ under this action.  The
cardinality of an orbit equals to the index of the
stabiliser.  It remains to note that the subgroup~$H_\phi$
is the stabiliser of the tail of $\phi$.

Let us prove the second assertion. Choose an element $x\in F$ of degree
one. A homomorphism $\alpha\:F\to G$ is uniquely determined by its tail
and the value $\alpha(x)$.  Moreover, for two homomorphisms $\alpha$ and
$\beta$ with the same tail, the quotient~$h=(\alpha(x))^{-1}\beta(x)$ must
stabilise this tail, i.e. $h$ must lie in $H_\alpha$. Indeed, for all
$f\in F$ of degree zero, we have
$$
\alpha(f^x)^h=
\alpha(f)^{\alpha(x)h}=
\alpha(f)^{\beta(x)}=\beta(f)^{\beta(x)}=\beta(f^x)=\alpha(f^x),
\qbox{i.e. $h$ centralises the subgroup $\alpha(\ker\deg)$;}
$$
and, for any element $f\in F$, we have
$$
\alpha(x)\alpha(f)H=
\alpha(xf)H=
\beta(xf)H=
\beta(x)\beta(f)H=
\alpha(x)h\beta(f)H=
\alpha(x)h\alpha(f)H,
\qbox{i.e. $h\in \alpha(f)H\alpha(f)^{-1}$.}
$$
Thus, $h=(\alpha(x))^{-1}\beta(x)\in H_\alpha$.

On the other hand, if $h$ is an arbitrary element of $H_\alpha$,
then the formula
$
f\mapsto\cases{
               \alpha(f),  &if $\deg f=0$\cr
               \alpha(x)h, &if $f=x$\cr
               }
$
defines a homomorphism with the same tail as $\alpha$ (by Lemma 0). This
homomorphism lies in $\Phi$ by Condition II of the Main Theorem.

Thus, for any $\alpha\in\Phi$, the set $\Phi$ contains precisely
$|H_\alpha|$ homomorphisms with the same tail as $\alpha$.  It
remains to note that, for similar homomorphisms $\psi$ and $\phi$, the
subgroups $H_\phi$ and $H_\psi$ have the same order, because they
are conjugate. This completes the proofs of assertion 2) and
the Main Theorem.

\REFERENCES


\[Bro00]
Brown K. S.
The coset poset and probabilistic zeta function of a finite group
//J. Algebra, 2000. V.225. P.989-1012.

\[Coll10]
Collins D. J.
Generating Sequences of Finite Groups.
Senior Thesis.
Cornell University Mathematics Department, 2010.
(Available here:
\newline
{\ttsmall
http://www.math.cornell.edu/m/sites/default/files/imported/Research/SeniorTheses/2010/collinsThesis.pdf}
)

\[Hall36]
Hall P.
The Eulerian functions of a group
// Quart. J. Math. Oxford Ser., 7 (1936), pp. 134-151.

\[HI\"O89]
Hawkes T., Isaacs I. M., \"Ozaydin M.
On the M\"obius function of a finite group
// Rocky Mountain J. Math. 1989. 19:4, 1003-1034

\[GRV12]
Gordon C., Rodriguez-Villegas F.
On the divisibility of $\#\Hom(\Gamma, G)$ by $|G|$
// J. Algebra. 2012. V.350, no.1, P. 300--307.
See also arXiv:1105.6066.

\[Iwa82]
S. Iwasaki,
A note on the $n$th roots ratio of a subgroup of a finite group
//
J. Algebra, 78:2 (1982), 460-474.

\[KM14]
Klyachko Ant. A, Mkrtchyan A. A.
How many tuples of group elements have a given property?
With an appendix by Dmitrii V. Trushin
// Intern. J. of Algebra and Comp., 2014, 24:4, 413-428.
See also arXiv:1205.2824

\[KT84]
Kratzer C., Th\'evenaz J.
Fonction de M\"obius d'un groupe fini et anneau de Burnside.
// Commentarii Mathematici Helvetici. 59:1(1984): 425-438.

\[Solo69]
Solomon L.
The solutions of equations in groups
// Arch. Math. 1969. V.20. no.3. P. 241--247.

\[Wils09]
Wilson R. A.
The Finite Simple Groups.
Graduate Texts in Mathematics.
Springer - 2009.

\endREFERENCES
\end